\documentclass[a4paper,11pt]{amsart}

\usepackage[english]{babel}
\usepackage{amsmath}
\usepackage{amsthm}
\usepackage{amssymb}
\usepackage{amscd}
\usepackage{hyperref}
\usepackage{geometry}
\usepackage{color}
\usepackage[foot]{amsaddr}

\title{Counting rational points and lower bounds for Galois orbits}
 
\date{}
\newtheorem{lemma}{Lemma}
\newtheorem{prop}{Proposition}

\newtheorem{thm}{Theorem}

\newtheorem{corollary}{Corollary}

\newtheorem*{conjecture}{Conjecture}

\newcommand{\Q}{\mathbb{Q}}
\newcommand{\C}{\mathbb{C}}
\newcommand{\R}{\mathbb{R}}

\newcommand{\G}{\mathbb{G}}
\author{Harry Schmidt}
\address{University of Manchester\\ harry.schmidt@manchester.ac.uk}
\begin{document}

\maketitle
\begin{abstract} In this article we present a new method to obtain polynomial lower bounds for Galois orbits of torsion points of one dimensional group varieties. 
 
 \end{abstract}
\section{Introduction}
In this article we introduce a method that allows in certain situations to obtain lower bounds for the degree of number fields associated to a sequence of polynomial equations. By a Galois bound for such a sequence we mean a lower bound for the degree of these number fields that is polynomial in the degree of the polynomials. 
Maybe the simplest example is 
\begin{align*}
X^n=1
\end{align*}
where $n$ runs through the positive integers. For a primitive solution $\zeta_n$ of this equation, that is, one that does not show up for smaller $n$ this is just the well-known Galois bound for roots of unity $[\Q(\zeta_n):\Q] \gg_\epsilon n^{1-\epsilon}$. The Galois bound here was proven by Gauss.

Another example are the equations 
\begin{align*}
B_n(X)=0
\end{align*}
 where $B_n$ is the division polynomial. For a fixed elliptic curve $E$ given in Weiertrass form 
\begin{align*}
E: Y^2=4X^3-g_2X-g_3
\end{align*}
this is defined by 
$[n](X,Y) = \left(\frac{A_n(X)}{B_n(X)},y_n\right)$ where $A_n \in \Q(g_2,g_3)[X]$ is a monic polynomial of degree $n^2$ and $B_n\in \Q(g_2,g_3)[X]$ is of degree $n^2-1$ with leading coefficient $n^2$ \cite[p.105, Exercise 3.7]{Silverman}. Here Galois bounds are not known as long as for roots of unity. The only known methods to obtain  Galois bounds so far seem to be either through Serre's open image theorem \cite[Corollary 9.5]{Lombardo} (in the non CM case), class field theory\cite[Th\'eor\`eme 2]{Remond}, \cite[(1.1)]{Silverberg}(in the CM case), or (in both cases) transcendence techniques as for instance applied  by  Masser \cite{Masser1} and further developed by David \cite{David}. 

For both of these examples we give a new proof of Galois bounds using essentially the same strategy for each (see Corollary \ref{galoisboundunity} and \ref{galoisboundelliptic}).

In joint work with Boxall and Jones we also consider fields obtained by adjoining all solutions $z$ of an equation $p^{\circ n}(z)=p^{\circ n}(y)$ for a fixed polynomial $p$ with coefficients in $\Q$ and certain fixed $y\in \Q$ (here $p^{\circ n}$ is $p$ iterated $n$ times)\cite{BJS}.  
\\

The use of counting of rational points on transcendental varieties in Diophantine geometry was first introduced by Pila and Zannier in \cite{PZ} to find yet another proof of Manin-Mumford and  initiated remarkable developments in Diophantine geometry. There are now some excellent accounts of these developments such as \cite{Scanlon} and \cite{Z} and we refrain from saying much more. \\

We only concentrate on the counting results. Here the basic idea is that for certain transcendental sets a subpolynomial bound for the number of rational points of height bounded by $H$ should hold.  The first result in this direction was proven by Bombieri and Pila \cite{BP}. Pila then later developed his determinant method to show among other things certain counting results for subanalytic surfaces \cite{Pila}. A conceptual jump was achieved with the introduction of $o$-minimal structures in the celebrated  Pila-Wilkie counting theorem \cite{PW}.\\

The question of improving the bound from subpolynomial to poly-logarithmic was perhaps around since the first types of such counting results were proven.  It was first shown by Surroca \cite{Surroca} that this is not possible in general. However, Wilkie conjectured that a poly-log bound should hold for the structure $\R_{\exp}$, the expansion of the reals by the real exponential function. 

\begin{conjecture}
Let $X^{trans}$ be the transcendental part of a  set $X\subset \R^k$ definable in $\R_{\exp}$. Then the following holds
\begin{align*}
\#\{x\in X^{trans}\cap \Q^k; H(x) \leq H\} \leq c_1\log H^{c_2},
\end{align*}
where $c_1,c_2$ are real constants depending on $X$.
\end{conjecture}
 We recall that $X^{trans}$ is $X$ deprived of all positive dimensional connected semi-algebraic sets contained in $X$.  
Now if we replace $\Q$ by a number field, the same type of bound would follow from the Conjecture. And we could ask how the constants depend on the number field.  Pila \cite{Ppfaff}, using essentially real analytic methods, proved a poly-log bound for Pfaffian curves and showed that they depend only on the degree of the number field. Jones and Thomas \cite{JT} then extended this to restricted Pfaffian surfaces. It is reasonable to conjecture that for Pfaffian varieties in any dimension the same should hold. In a landmark work \cite{BN} Binyamini and Novikov, who combined complex analytic methods with Khovanskii's zero-estimates,  proved a poly-log bound for sets definable by restricted elementary functions. Another method to prove poly-log bounds for analytic functions using Siegel's lemma  was introduced by Masser \cite{zeta}. And a variant due to Wilkie  \cite{WilkieSiegel} was used to give an alternative proof of his theorem with Pila. See also Habegger's work \cite{Habeggerapproximation} where this approach is used in place of the determinant method, in order to count algebraic points near definable sets.  Masser's result on the Riemann $\zeta$ function seems to be the first with a polynomial dependence on the degree of the number field \cite[p.2045 (15)]{zeta}. Related results for analytic functions were obtained by  Boxall and Jones \cite{BJ2}, \cite{BJ1}, Besson \cite{Besson} and Jones and Thomas \cite{JTzeta}. 
A reasonable question is whether the dependence on the degree of the number field can be made polynomial in the structures we mentioned. So for $X$ as above one might ask whether the following holds
\begin{align*}
\#\{x\in X^{trans}\cap K^k; H(x) \leq H\} \leq c_1 d^{c_3}\log H^{c_2},
\end{align*}
where $d=[K:\Q]$ and $c_1,c_2,c_3$ depend on $X$. In fact we would expect that a bound as above holds for $X$ definable by an extension of the reals by the complex exponential function restricted to a certain fundamental domain or the restricted $j$-function. Even more generally the above should hold for $X$ definable in the extension of the reals by the restriction of a uniformization map of a mixed Shimura variety on a suitable fundamental domain. 

It seems that the methods leading to  a poly-log bound always lead to a polynomial dependence on the degree. We will demonstrate this for Pfaffian curves and restricted Pfaffian surfaces. In what follows we write $X(K,H)$ for the points of $X$ with coordinates in the number field $K$ and multiplicative Weil  height at most $H$.  (The height is extended to tuples by taking the maximum).  \\

Pila introduced the notion of mild parametrization. A set has a $(J,A,C)$-mild parametrization if it can be covered by the image of $J$, $C^{\infty}$ functions on the unit cube whose partial derivatives of order $\mu$ are bounded by $\mu!(A|\mu|^C)^{|\mu|}$  \cite[Definition 2.1 and 2.4]{Pmild}.  Pila showed that the rational points of a  set $X$ that has mild parametrization can be covered by poly-log many hypersurfaces. In his work one can make the dependence on the number field polynomial.
\begin{thm}\label{mild} Let $X \subset (0,1)^n$ have a $(J,A,C)$-mild parametrization. There exist effectively computable constants  $C_1,C_2,C_3$ depending only on $J,A,C,n$ and $\dim X$ such that $X(K,T)$ is contained in the union of $C_1d^{C_2}\log T^{C_3}$ hypersurfaces of degree bounded by 
$(d^2\log T)^{\dim X/(n-\dim X)}$.
\end{thm}
Pila's approach was then used by Jones and Thomas to prove a poly-log bound for any transcendental implicitly defined Pfaffian curve. In their work \cite{JT} the dependence on the degree of the number field can again be made polynomial.  First for curves.
\begin{thm} \label{curve}
Suppose that $I$ is an open interval in $\R$ and that $f:I \rightarrow \R$ is a transcendental implicitly defined Pfaffian function of complexity $(n,r,\alpha, \beta)$. Then for $T\geq e$, and the graph $X$ of $f$
\begin{align*}
|X(K,H)|\leq  c(n,r,\alpha,\beta)d^{6n +6r +15}\log H^{3n +3r +8}.
\end{align*}
where $c(n,r,\alpha,\beta)=2^{r(r-1)}6^{\frac52n+r + \frac{53}2}(n+2)^{n+3r+1}(\alpha+\beta)^{2n+2r+1}$.
\end{thm}

They used this result and a stratification result due to Khovanskii and Vorobjov to show a poly-log bound for restricted Pfaffian surfaces. Again one can get a polynomial dependence on the degree from their work. 
First let $X \subset (0,1)^n$ be a restricted semi-Pfaffian surface and assume that $X$ has a mild parametrization with parameters $J,A,C$ (with $J,A,C$ as above) bounded by a constant $M$ and that $X$ is not contained in an algebraic hyper-surface.   From Jones and Thomas work \cite{JT} one can deduce the following.
\begin{thm}\label{surface} There are effectively computable constants $C_1,C_2,C_3$ depending only on the format of $X$ and on $M$ such that
\begin{align*}
|X^{trans}(K,T)|\leq C_1d^{C_2}\log T^{C_3}.
\end{align*}
\end{thm}

The assumption on $X$ to not be contained in an algebraic hypersurface could be dropped but we will omit the technical details.\\

Here is how this article is structured. In the next section \ref{pfaff} we show how to deduce Theorem \ref{curve} and \ref{surface} from the work of Pila and Jones and Thomas. Then in section \ref{torsion} we apply these results to find a new proof of polynomial lower Galois bounds for torsion points of $\G_m$ and elliptic curves. \\

The author would like to thank Philipp Habegger, Gareth Jones and Jonathan Pila for many encouraging discussions. He would also like to thank Philipp Habegger, Gareth Jones and David Masser for pointing out several typos in a previous draft and for their comments leading to an improvement in the exposition. He would also like to heartily thank the referee for his comments. Finally  he also thanks the Engineering and Physical Sciences Research Council for support under grant  EP/N007956/1. 

\section{Pfaffian curves and surfaces}\label{pfaff}
We start by recalling work of Pila on Pfaffian curves. We will use the notions from Definition 6.4 in \cite{Pselecta}. 
 We also define $X^{size}(K,H)$ to be the set of points in $X$ with coordinates in  $K$-points such that $H^{size}(x) \leq H$ where $H^{size}(\alpha)$ for an algebraic number $\alpha$ is defined as
$H^{size}(\alpha) =\max_\sigma\{den(\alpha), |\sigma(\alpha)|\}$
 where we take the maximum over all embeddings $\sigma$ of $\Q(\alpha)$ into $\C$ and $den(\alpha)$ is the smallest positive integer $\gamma$ such  $\gamma  \alpha$ is an algebraic integer (see \cite[Definition 6.3]{Pselecta}). We extend $H^{size}$ to tuples by taking the maximum.   
\begin{lemma}\label{lemma1} Let $K$ be a real number field with $[K:\Q]=d$. Let $\delta\geq 1, T \geq 1, L \geq 1/T^{4d}$. Put $D = (\delta+1)(\delta+2)/2$. Let $I$ be an interval of length at most $L$ and $f:I \rightarrow \R$ a function. Suppose that $f$ has $D-1$ continuous derivatives on $I$, with $|f'| \leq 1$, and let $X$ be the graph of $f$ on $I$. Then $X^{size}(K,T)$ is contained in the union of at most 
\begin{align*}
6(D!)^{\frac{2d}{D(D-1)}}(LT^{4d})^{\frac{4}{3(\delta+3)}}A_{L,D-1}(f)
\end{align*}
plane algebraic curves of degree at most $\delta$. 
\end{lemma}
\begin{proof} The proof is the same as the proof of \cite[Lemma 6.5]{Pselecta} except that instead of estimating   $(D!)^{2/(D(D-1))}$  we leave it as it is and estimate $D^{4d/(3(D-1))}\leq 5$ there . 
\end{proof}
\begin{prop}\label{implicitp} Let $\delta\geq 1, D=(\delta+1)(\delta+2)/2, d \geq 1, T\geq e, L \geq 1/T^{4d}$ and $I\subset \R$ an interval of length $\leq L$. Let $K\subset \R$ be a number field of degree $d$. Let $f:I \rightarrow \R$ have $D$ continuous derivatives, with $|f'|\leq 1$ and $f^{(j)}$ either non-vanishing in the interior of $I$ or identically vanishing, for $j=1, \dots, D$. Let $X$ be the graph of $f$. Then $X^{size}(K,T)$ is contained in the union of at most 
\begin{align*}
54D(D!)^{2d/(D(D-1))}(LT^{4d})^{4/(3(\delta+3))}\log (eLT^{4d}).
\end{align*}
real algebraic curves of degree at most $\delta$. 
\end{prop}
\begin{proof} We  follow the proof of  \cite[Prop. 6.7]{Pselecta} word by word but using the estimate from Lemma \ref{lemma1} instead of the one given in Lemma 6.5 there. 
\end{proof}

Finally we state a result of Jones and Thomas a little more explicitly. Recall their definition of implicitly defined Pfaffian function on  page 640 of \cite{JT} with its notion of complexity. 
\begin{thm} \label{implicit}
Suppose that $I$ is an open interval in $\R$ and that $f:I \rightarrow \R$ is a transcendental implicitly defined Pfaffian function of complexity $(n,r,\alpha, \beta)$. Then for $T\geq e$, and the graph $X$ of $f$
\begin{align*}
|X^{size}(K,T)|\leq  2^{r(r-1)}6^{\frac52n+r + \frac{53}2}(n+2)^{n+3r+1}(\alpha+\beta)^{2n+2r+1}d^{3n +3r +8}\log T^{3n +3r +8}.
\end{align*}
\end{thm}
\begin{proof} We follow their  proof and for $c_3,c_4$ there it is not hard to estimate
\begin{align*}
c_3c_4d^{n+r+1}D^{n+r+2}\leq 2^{r(r-1)}6^{\frac32(n+3) + 3}(n+2)^{n+3r+1}(\alpha+\beta)^{2n+2r+1}d^{n+r+1}D^{n+r+2}\\
\leq 2^{r(r-1)}6^{\frac52n+r + \frac{23}2}(n+2)^{n+3r+1}(\alpha+\beta)^{2n+2r+1}d^{3n+3r+5}.
\end{align*}
 Now setting $\delta= [d\log T]$ and $L=2T$ in Proposition \ref{implicitp} we get that the number of hypersurfaces can be bounded by 
\begin{align*}
6^{15}d^3\log T^3
\end{align*}
(note the slightly confusing fact that $d$ is now the degree of the number field) and we obtain the theorem.
\end{proof}
Now we deduce Theorem \ref{curve} from Theorem \ref{implicit} by noting that $X(K,T)\subset X^{size}(K,T^d)$.

We continue by recording Pila's work \cite{Pila}. Recall the definition of $(J,A,C)$ mild \cite[Definition 2.1 and 2.4]{Pmild}.

Pila proved that for $X$ admitting a mild parametrization the set $X(K,T)$ is contained in a union of hypersurfaces whose number and degree we can control in terms of $J,A,C$. \\

In order to get a polynomial dependence on the degree of $K$  we we are going to make a slight adjustment at one of the steps in Pila's proof (which was more concerned with the dependence on $T$) and show that we can extract a polynomial dependence on the degree. 
\begin{proof}[\textit{Proof of Theorem \ref{mild}} ]We follow the proof in \cite{Pmild} right up to the choice of the degree of the hypersurface on p.503. The degree of the number field there is $f$. Instead of choosing $d$ there equal to $[\log T^{k/(n-k)}]$ we choose it to be equal to $[(f\log T)^{k/(n-k)}]$ there which kills the $f$ in the exponent.  Then in the corollary there we get  $f^2$.
\end{proof}

We start with our investigation of surfaces. For this we use the notion of a semi-Pfaffian set as in \cite[p.640]{JT} with its notion of format.    First we record that the following holds. Let $X$ be a connected semi-pfaffian surface in $(0,1)^n$ with a mild parametrization with parameters bounded by a constant $M_X$ that is not contained in any algebraic hypersurface. 
\begin{prop}\label{polysurface} There exist effectively computable constants $c_1,c_2,c_3,c_4$ depending only on the format of $X$ and $M_X$ such that for any hypersurface $Z$ in $\R^n$ of degree $d_Z$ holds
\begin{align*}
|(X\cap Z)^{trans}(K,T)|\leq c_1d_Z^{c_2}d^{c_3}\log T^{c_4}
\end{align*}
\end{prop}
\begin{proof} We first note that each component of $X\cap Z$ has dimension at most 1 since $X$ is not contained in $Z$. Furthermore, by Khovanskii's zero-estimates the number of connected components of the intersection $X\cap Z$ grows polynomially in $d_Z$ with the growth depending only on the format of $X$ as is already pointed in the proof of Proposition 5.3 of \cite{JT}.   Now if the component is a point the counting becomes straightforward. So we may assume it is a curve. If the curve is algebraic it does not belong to $(X\cap Z)^{trans}$. If the curve is transcendental there is some projection to $\R^{3}$ that is a transcendental curve as well. For each projection we  can follow the proof of \cite[Proposition 5.3]{JT} line by line but in the displayed equation just after (7) we use Theorem \ref{implicit} to get a polynomial dependence on $d$. There are $n\choose{3}$ such projections so we need to multiply the final estimate by this number as well.  
\end{proof}

Combining Theorem \ref{implicit} with Proposition \ref{polysurface} we obtain Theorem \ref{surface}.

\section{Torsion points}\label{torsion}

\begin{corollary}\label{galoisboundunity} For $n\geq 4$ holds 
\begin{align*}
[\Q(\zeta_n):\Q]\geq n^{1/40}\log n^{-1/2}/6.
\end{align*}
\end{corollary}
\begin{proof}We consider the function $\cos(2\pi \theta)$ on the interval $(-\frac12,\frac12)$. This is Pfaffian of degree $(2,1)$ and order 2. For $n\geq 4$ let $\zeta_n$ be a primitive $n$-th root of unity. Then $\Re(\zeta_n) =\cos(2\pi (k/n))$ for some integer $k$ in $(-n/2,n/2)$ and $(k,n)=1$. Since $\overline{\zeta_n}= 1/\zeta_n$ this lies in $\Q(\zeta_n)$. We also have that $\zeta_n^l \in \Q(\zeta_n)$ for an integer $1\leq l\leq n-1$ and that $\Re (\zeta_n^l) = \cos(2\pi l'/n)$ for some $l'\in (-n/2,n/2)$ if we exclude $l =n/2$. From elementary height inequalities follows $H(\cos(2\pi k/n))\leq 4$. We set $d=[\Q(\zeta_n):\Q]$ and conclude from Theorem \ref{implicit} that 
\begin{align*}
n-2\leq c(2,2,2,1)d^{40}\log n^{20}
\end{align*}
where $c(2,2,2,1)\leq 6^{50}$.
\end{proof}
Of course this is far from the real lower bound but the proof uses the same method that we will use to deduce lower bounds for torsion points of elliptic curves.

In order to keep the exposition short we refrain from proving more explicit bounds. We will only show that the dependence on the elliptic curve is only on the height of the elliptic curve. \\

We fix a lattice $\Lambda\subset \C$ with generators $\omega_1,\omega_2$.  Let $\wp_\Lambda$ be the Weierstrass function associated to $\Lambda$ and $E$ the associated elliptic curve. In what follows we denote complex conjugation by an upper bar.  We define the function $f_\Lambda$ by 
\begin{align*}
f_\Lambda(b_1,b_2) =(\Re{\wp_\Lambda(b_1\omega_1+b_2\omega_2)},\Im\wp_\Lambda(b_1\omega_1+b_2\omega_2));~~ (b_1,b_2)\in [0,1)^2, b_1^2+b_2^2\neq 0.
\end{align*}
Let $X_\Lambda$ in $\R^4$ be the graph of $f_\Lambda$. For any algebraic hypersurface $Z$ in $\R^4$ the following holds.
\begin{lemma} \label{transcendence} Any positive dimensional component of $Z\cap X_\Lambda$ is a  transcendental curve. 
\end{lemma}
\begin{proof}  It is enough to prove that the intersection $Z\cap X_\Lambda$ has no 2 dimensional components and that the one dimensional components are transcendental curves. Suppose first that there exists a 2 dimensional component $U$. Then there exist complex analytic functions $r_1,r_2$ in some poly-disc in $\C^2$ such that the determinant of the Jacobian of $(r_1,r_2)$ does not vanish and such that their restriction to $\R$ is real. Further $trdeg_\C \C(r_1,r_2, \wp_\Lambda(r_1\omega_1+r_2\omega_2), \wp_{\overline{\Lambda}}(r_1\overline{\omega}_1+r_2\overline{\omega}_2))\leq 3$. By Ax's theorem \cite{Ax} this implies that $y=(\wp_\Lambda(r_1\omega_1+r_2\omega_2),\wp_{\overline{\Lambda}}(r_1\overline{\omega}_1+r_2\overline{\omega}_2))$ is contained in a translate of an algebraic  subgroup of $E\times\overline{E}$. Thus there exists an isogeny $\beta:\overline{E}\rightarrow E$ and endomorphisms $\alpha_1, \alpha_2$ that act by multiplication by a non-zero complex number on the tangent space that we again denote by $\beta, \alpha_1, \alpha_2$ such that $\alpha_1(r_1\omega_1 + r_2\omega_2) + \alpha_2 \beta(r_1\overline{\omega}_1+ r_2\overline{\omega}_2) = 0 \mod \C $. Further $\beta$  acts by sending $b_1\overline{\omega}_1+b_2\overline{\omega}_2$ to $b_1^B\omega_1+b_2^B\omega_2$ where $(b_1^B,b_2^B)= (b_1,b_2)B$ with $B$  an integer matrix. If we set $(a,b) = (1,0)B, (c,d) = (0,1)B$ then $\overline{\tau} = \frac{a\tau + b}{c\tau + d}$ for $\tau = \omega_1/\omega_2$ and so $\Im(\overline{\tau}) = \frac{\det B\Im(\tau)}{|c\tau + d|^2}$. Thus $B$ has negative determinant.  There is thus a relation of the form 
\begin{align*}
(r_1,r_2)A_1+ (r_1,r_2)BA_2 =0 \mod \R^2.
\end{align*}
where $A_1,A_2$ are integer matrices given by the action of $\alpha_1, \alpha_2$. These have positive determinant by an argument as above for $B$. It follows that the matrix $A_1+BA_2$ has rank at least 1 and so $r_1,r_2$ are linearly related over $\Q \mod \R$ which contradicts our assumption that the Jacobian of $(r_1,r_2)$ is non-singular. Thus there are no 2-dimensional components. 

Now let $U$ be a 1 dimensional component and assume that it is algebraic. Thus there are complex analytic functions $r_1,r_2$ not both constant such that $trdeg_\C \C(r_1,r_2, \wp_\Lambda(r_1\omega_1+r_2\omega_2), \wp_{\overline{\Lambda}}(r_1\overline{\omega}_1+r_2\overline{\omega}_2))\leq 1$. Now as above this implies that $r_1,r_2$ are linearly related over $\Q$, $\mod \R$. That is we may assume that we can write $r_2=c_1r_1 +c_2$ for $c_1,c_2\in \R$ and that $r_1$ is not constant. Setting $z_r= r_1(\omega_1+c_1\omega_2) +c_2\omega_2$ we have $trdeg_\C\C(z_r, \wp_\Lambda(z_r))\leq 1$. By Ax's theorem this implies that $z_r$ is constant 
and so $\omega_1+c_1\omega_2=0$ which is absurd since $\Im(\omega_2/\omega_1)\neq 0$.  This proves the claim. 
\end{proof}


We borrow some estimates from  Masser-W\"ustholz. First \cite[Lemma 3.2]{MW}.  
\begin{lemma} \label{upperbound} There exists an effectively computable absolute constant $C$ such that 
\begin{align*}
|\wp_\Lambda(z) -\wp_\Lambda(\omega_2/2)| \leq Cd(z,\Lambda)^{-2}.
\end{align*}
where $d(z,\Lambda)$ is the minimal distance of $z$ to an element of $\Lambda$. 
\end{lemma}
 
Pick generators $\omega_1, \omega_2$ of $\Lambda$ such that $\tau =\omega_2/\omega_1$ satisfies $|\Re(\tau)|\leq \frac12, |\tau|\geq 1$. 
\begin{lemma}\label{lowerbound} For any $\epsilon >0$ let $B_\epsilon$ be the box consisting of $z=t_1\omega_1 +t_2\omega_2$ with $|t_1|\leq \frac12, |t_2|\leq 1/2-\epsilon$. There is an effectively computable constant $C_\epsilon$ depending only on $\epsilon$ such that 
\begin{align*}
|1/(\wp_\Lambda(z)-\wp_\Lambda(\omega_2/2))| \leq C_\epsilon\exp(\pi \Im(\tau)), z \in B_\epsilon. 
\end{align*}
\end{lemma}

\begin{proof} We consider the expansion \cite[(3.3)]{MW}. From the proof there follows that the absolute value of $F(1/2)\prod_{n=1}^{\infty}\{ F(n)/F(n-\frac12)\}$ is bounded from above. For $F(1/2)$ we note that for $q,Q$ there 
\begin{align*}
|1-q^{\frac12}Q^{\pm2}|\geq 1- \exp(-2(1/2 \pm t_2)\pi\Im(\tau)) \gg_\epsilon\exp(-2\epsilon \pi\Im(\tau)).
\end{align*}
The last thing to check is that $|\sin w|\ll\exp (2\pi(1/2- \epsilon) \Im(\tau))$.
\end{proof}

We pass to the Legendre family $E_\lambda$ and set $X_\lambda = \wp_\lambda +\frac13(\lambda+1)$ where $\wp_\lambda$ is associated to the lattice $\Lambda_\lambda$ generated by the differential $\frac{dX}{2\sqrt{X(X-1)(X-\lambda)}}$. 
Now assume that $\lambda$ satisfies 
\begin{align}\label{condition}
|\lambda|\leq 1, |1-\lambda|\leq1, \Re(\lambda)\leq \frac12
\end{align}
 and set $\omega_1, \omega_2$ to be given by hypergeometric series such as in \cite[p.5]{pfaffian}. It can be checked that $X_\lambda(z) = \wp_\lambda(z) -\wp_\lambda(\omega_2/2)$ (see for example \cite[Lemma 5.1]{Schmidtadditive}). 

We first define $U_1$ be given by 
\begin{align*}
U_1= \{z =b_1\omega_1+b_2 \omega_2;~~ b_2 \in [1/30,29/30] ,b_1\in [0,1] \}
\end{align*}
 \begin{lemma} \label{distance} We have
\begin{align*}
d(z,\Lambda_\lambda) \geq |\omega_2|/60,~~z \in U_1.
\end{align*}
\end{lemma}
\begin{proof}  For $z = b_1\omega_1 + b_2\omega_2 \in U_1$ we have that for $\omega \in \Lambda_\lambda$, $z - \omega = c_1 \omega_1 + c_2\omega_2$ for $c_1,c_2 \in \R$ and $|c_2| \geq 1/30  $. Now $z-\omega = (c_1 + c_2\tau)\omega_1$ for $\tau = \omega_2/\omega_1$ and 
\begin{align*}
|z-\omega| = |\omega_1||c_1 +c_2\tau| \geq |\omega_1||c_2||\Im(\tau)|\geq |\omega_2|/60.
\end{align*}
\end{proof}
We define $U_2= \{z=b_1\omega_1 + b_2\omega_2; b_1,b_2\in [-29/60,29/60]\}$. \\

After a calculation one finds that for any $\lambda' \in \C\setminus \{0,1\}$ we can find $\lambda$ satisfying (\ref{condition}) such that  one of the following holds 
\begin{align}\label{lattice}
\Lambda_{\lambda'} = \epsilon \Lambda_\lambda, ~~ \Lambda_{\lambda'} = \epsilon (1-\lambda)^{\frac12}\Lambda_{\lambda}, ~~ \Lambda_{\lambda'} = \epsilon \lambda^{\frac12}\Lambda_{\lambda},
\end{align}  
where $\epsilon \in \{1,i\}$. 
This follows from the fact that the transformation $\lambda \rightarrow 1-\lambda$ scales the lattice by $i$ while $\lambda \rightarrow 1/\lambda$ scales the lattice by a factor of $\sqrt{\lambda}$ \cite[(5),(8)]{Fettis}.  
Since the generators $\omega_1, \omega_2$  are  such that $\omega_2/\omega_1 = \tau$ lies in the standard fundamental domain as above \cite[Lemma8]{pfaffian} the same holds for $ \epsilon \lambda^{\frac12}\omega_1, \epsilon \lambda^\frac12\omega_2$ respectively $\epsilon (1-\lambda)^\frac12\omega_1, \epsilon (1-\lambda)^\frac12 \omega_2$. And for each $\lambda' \in \C\setminus\{0,1\}$ we pick generators $\omega_1',\omega_2'$ to be equal to such a pair.\\
 
If we can choose $\lambda$ in (\ref{condition}) such that $\Lambda_{\lambda'} = \epsilon \Lambda_{\lambda}$  or $\Lambda_{\lambda'} = \epsilon (1- \lambda)^{\frac12}\Lambda_{\lambda}$ we set 
\begin{align*}
f_{\lambda'}(b_1,b_2)  = X_\lambda(b_1\omega_1 +b_2\omega_2), z=b_1\omega_1 +b_2\omega_2 \in U_1\cup \{r\omega_1; r \in (0,1)\}.
\end{align*}
Otherwise  we set 
\begin{align*}
f_{\lambda'}(b_1,b_2) & = \lambda/X_{\lambda}(b_1\omega_1 +b_2\omega_2),  z=b_1\omega_1 +b_2\omega_2 \in U_2.
\end{align*}

We define $X_{\lambda'}(z)  = \wp_{\lambda'}(z) - \wp_{\lambda'}(\omega_2'/2)$ and from the homogeneity of the Weierstrass $\wp$ function follows that $X_{\lambda'}(b_1\omega_1' + b_2 \omega_2') = \pm f_{\lambda'}$ if $\Lambda_{\lambda'} = \epsilon \Lambda_{\lambda}$ while  if $\Lambda_{\lambda'} = \epsilon (1-\lambda)^{\frac 12} \Lambda_{\lambda}$ then $X_{\lambda'} = \pm f_{\lambda'}/(1-\lambda)$ and finally if $\Lambda_{\lambda'}  = \epsilon \lambda^{\frac12}\Lambda_\lambda$ then $X_{\lambda'}(b_1\omega_1' + b_2\omega_2')= \pm 1/f_{\lambda'}$.

\begin{lemma} \label{mildpar}There exists an integer $T$ and  effectively computable absolute constants $A_1, A_2$ (not depending on $\lambda'$) such that the graph $G_{\lambda'}$ of $f_{\lambda'}/T$ restricted to $(b_1,b_2)$ such that $z\in U_1$ respectively $z\in U_2$ has a $(1,A_1,A_2)$ mild parametrization. 
\end{lemma}
\begin{proof} 
First suppose that $\Lambda_{\lambda'}= \epsilon \Lambda_{\lambda}$ or $\Lambda_{\lambda'}= \epsilon(1-\lambda)^\frac12 \Lambda_{\lambda}$. We use Cauchy's formula for the $n$-th derivative ($n \geq 1$)
\begin{align*}
X_{\lambda}^{(n)}(z)/n! =\frac1{2\pi i} \oint\frac{X_\lambda(w)dw}{(z-w)^{n+1}} 
\end{align*}
and integrate along the circle $|z-w| =d(z,\Lambda_\lambda)/2$. By Lemma \ref{distance} $d(z,\Lambda_\lambda)/2\geq |\omega_2|/120$ for $z\in U_1$. So by Lemma \ref{upperbound} there exist absolute constants $\tilde{A}_1,\tilde{A}_2$ such that 
\begin{align*}
|\omega_2|^n|X_{\lambda}^{(n)}| \leq \tilde{A}_1(\tilde{A}_2)^{n}n!, \text{ for } z\in U_1.
\end{align*}
As the absolute value of $\omega_1$ is bounded by an absolute constant (see for example \cite[Lemma 12]{pfaffian} ),  if we pick an integer $T$ whose absolute value is greater than the maximum absolute value of $X_\lambda$ on $U_1$ we find that a mild parametrization of $G_{\lambda'}$ is given by 
\begin{align*}(t_1,t_2)\rightarrow (t_1, 1/30 + (28/30)t_2, \Re f_{\lambda'}(t_1,1/30 + (28/30)t_2)/T,\Im f_{\lambda'}(t_1,1/30 + (28/30)t_2)/T). 
\end{align*}

Now assume that $\Lambda_{\lambda'}= \epsilon\sqrt{\lambda}\Lambda_\lambda$. We first note that from the Fourier expansion of $\lambda$ follows that $|\lambda|\ll \exp(-\pi \Im(\tau))$ \cite[p.117]{Chandrasekharan} so we find with Lemma \ref{lowerbound} that $\lambda/X_\lambda$ is bounded by  an aboslute constant on $B_{1/120}$. We again use Cauchy's formula but this time we first note that using the same arguments as in Lemma \ref{distance} the minimal distance of an element in $U_2$ to the boundary of $B_{1/120}$  is bounded from below by an absolute constant $c_1$. 
We can then pick the circle $|z-w| =c_1/2$ and using Cauchy's formula find that  we can pick $\tilde{A}_1, \tilde{A}_2$ such that 
\begin{align*}
|(\lambda/X_\lambda)^{(n)}(z) | \leq \tilde{A}_1\tilde{A}_2^nn!.
\end{align*}
Now it remains to note again that $\omega_1$ is bounded absolutely while $|\omega_2| \ll -\log |\lambda|$ and so $|\lambda|^{\frac12}|\omega_2|$ is also bounded by an absolute constant. We take  $T$ to be also larger than the maximum of $|\lambda/X_\lambda|$ on $U_2$. As in the previous case  we have established the mild parametrization of $G_{\lambda'}$. 
\end{proof}

Now let $\Gamma_{\lambda'}$ be the graph of $f_{\lambda'}$.  
\begin{corollary}\label{Weierstrasscounting} There exist effectively computable  absolute constants $\gamma_1, \gamma_2, \gamma_3$ such that 
\begin{align*}
|\Gamma_{\lambda'}(K,H)| \leq \gamma_1d^{\gamma_2}\log H^{\gamma_3}.
\end{align*}
\end{corollary}
\begin{proof} Theorem 1 of \cite{pfaffian} implies that $\Gamma_{\lambda'}$ is a finite union of semi-Pfaffian surfaces with the entries of its format and the number of surfaces in the union bounded by a constant independent of $\lambda'$.  For $\lambda'$ such that $\Lambda_{\lambda'}=\epsilon \sqrt{\lambda}\Lambda_{\lambda}$ the corollary follows directly (after rescaling) from Lemma \ref{mildpar}, Lemma \ref{transcendence} and Theorem \ref{surface}. For $\lambda'$ such that $\Lambda_{\lambda'}=\epsilon \Lambda_{\lambda}$ or $\Lambda_{\lambda'}= \epsilon (1-\lambda)^{\frac12}\Lambda_{\lambda}$ we need to also infer the use of Theorem \ref{curve} for the piece given by $f_{\lambda'}$ restricted to $(0,1)$. 
\end{proof}
We note here that there is a certain uniformity in the counting, since the constants $\gamma_1, \gamma_2, \gamma_3$ do not depend on $\lambda$.\\

From the corollary follows another corollary.
\begin{corollary}\label{galoisboundelliptic}  For algebraic $\lambda'$, let $P \in E_{\lambda'}(\overline{\Q})$ be a torsion point of order $n$ and $d=[\Q(P,\lambda'):\Q]$. There exist effectively computable positive  absolute constants $\delta_1, \delta_2, \delta_3$ such that 
\begin{align*}
d\geq \delta_1(1+h(\lambda'))^{-\delta_2}n^{\delta_3}. 
\end{align*}
\end{corollary}
\begin{proof} 
For each $\lambda' \in \C\setminus\{0,1\}$ there is a $T \in \{0,1,\lambda\}$ such that the abscissa of $P$ is equal to $X_{\lambda'}(z) - T$ where $z$ is a logarithm of $P$ (with respect to the tangent space at the identity). If $P$ is torsion of order $n$ then we can pick $z  = \frac kn \omega_1 + \frac ln\omega_2$ where $k,l$ are positive integers not larger than $n$ and such that $(k,l,n) =1$.  Clearly $X_{\lambda'}(z)$ lies in the field $\Q(\lambda',P)$ and the same holds for $X_{\lambda'}(kz), k = 1,\dots, n-1$. 
The logarithmic Weil height of the abscissa of a torsion point is bounded by $c(1+h(\lambda'))$ (\cite[p.40, Theorem]{Zimmer} or \cite[p.467]{six}) for $c$ absolute. An elementary computation shows that the set consisting of $m(\frac kn\omega_1 + \frac ln\omega_2) \mod \Lambda_{\lambda'}, m=1, \dots,n-1$ has $(1/c)n$ representatives in $U_1\cup\{ r\omega_1;r\in (0,1)\}$ respectively $U_2$. The difference between the heights of $f_{\lambda'}$ and $X_{\lambda'}$ is bounded absolutely and so  we find $(1/c)n$ points in $\Gamma_{\lambda'}(\Q(\lambda', P),H)$ with $\log H\leq c(1+h(\lambda'))\log n$.  From Corollary \ref{Weierstrasscounting} follows  the present corollary.
\end{proof}

\bibliography{Galoisbounds}
\bibliographystyle{amsalpha}

\end{document}